\documentclass[a4paper,11pt]{article}
\usepackage{amsmath,amssymb,amscd}
\def\qedsymbol{q.e.d.}
\def\proofname{Proof.}
\newenvironment{Proof}{\par\noindent{\it\proofname}}{{\unskip\nobreak\hfill{\it\qedsymbol}}\par\vskip 9pt}
\newenvironment{Proof*}{\par\noindent}{{\unskip\nobreak\hfill{\it\qedsymbol}}\par\vskip 9pt}
\ifx\undefined\bysame\newcommand{\bysame}{\leavevmode\hbox to3em{\hrulefill}\,}\fi
\numberwithin{equation}{section}
\newtheorem{Thm}{Theorem}[section]
\newtheorem{Lemma}[Thm]{Lemma}
\newtheorem{Def}[Thm]{Definition}
\newtheorem{Prop}[Thm]{Proposition}

\newtheorem{Cor}[Thm]{Corollary}

\newtheorem{Remark}[Thm]{Remark}

\def\bff{\boldsymbol F}
\def\ca{\mathcal A}
\def\vzr{\boldsymbol 0}
\def\alg{{\mathcal A}{\it lg}}
\def\Part{\mathrm{Part}}
\def\set{{\mathcal S}\!{\it et}}
\def\seq{\mathrm{Seq}}
\def\spec{\mathrm{Spec}}
\def\sq{{\textsl Sq}}

\begin{document}

\title{The Steenrod algebra from the group theoretical viewpoint}

\author{Atsushi Yamaguchi}

\date{}

\maketitle

\begin{abstract}
In the paper ``The Steenrod algebra and its dual'' \cite{[17]}, J.Milnor determined the structure of 
the dual Steenrod algebra which is a graded commutative Hopf algebra of finite type. 
We consider the affine group scheme $G_p$ represented by the dual Hopf algebra of the mod $p$ 
Steenrod algebra. 
Then, $G_p$ assigns a graded commutative algebra $A_*$ over a prime field of finite characteristic 
$p$ to a set of isomorphisms of the additive formal group law over $A_*$, 
whose group structure is given by the composition of formal power series (\cite{ay1} Appendix). 
The aim of this paper is to show some group theoretic properties of $G_p$ by making use of this 
presentation of $G_p(A_*)$. 
We give a decreasing filtration of subgroup schemes of $G_p$ which we use for estimating the length 
of the lower central series of finite subgroup schemes of $G_p$. 
We also give a successive quotient maps $G_p\xrightarrow{\rho_0}G_p^{\langle1\rangle}
\xrightarrow{\rho_1}G_p^{\langle2\rangle}\xrightarrow{\rho_2}\cdots\xrightarrow{\rho_{k-1}}
G_p^{\langle k\rangle}\xrightarrow{\rho_k}G_p^{\langle k+1\rangle}\xrightarrow{\rho_{k+1}}\cdots$ of 
affine group schemes over a prime field $\bff_p$ such that the kernel of $\rho_k$ is a maximal 
abelian subgroup. 
\end{abstract}

\section{Introduction}

The ring of cohomology operations on the mod $p$ ordinary cohomology theory which is denoted by 
$\ca_p^*$ below is called the Steenrod algebra. 
We denote by $\ca_{p*}$ the dual Hopf algebra of $\ca_p^*$ whose structure is determined by 
J. Milnor in \cite{[17]} as follows. 
$$\ca_{p*}=\begin{cases}
\bff_2[\zeta_1,\zeta_2,\dots,\zeta_i,\dots] &\text{if $p=2$}\\
E(\tau_0,\tau_1,\dots,\tau_i,\dots)\otimes_{\bff_p}\bff_p[\xi_1,\xi_2,\dots,\xi_i,\dots] &
\text{if $p$ is an odd prime}\end{cases}$$
Here, we set $\xi_0=\zeta_0=1$ and $\deg\,\zeta_i=2^i-1$, $\deg\,\tau_i=2p^i-1$, 
$\deg\,\xi_i=2(p^i-1)$. 
The coproduct $\mu:\ca_{p*}\to\ca_{p*}\otimes_{\bff_p}\ca_{p*}$ is given by the following formulas. 
$$\mu(\zeta_n)=\sum\limits_{k=0}^n\zeta_{n-k}^{2^k}\otimes\zeta_k, \qquad 
\mu(\xi_n)=\sum\limits_{k=0}^n\xi_{n-k}^{p^k}\otimes\xi_k, \qquad
\mu(\tau_n)=\sum\limits_{k=0}^n\xi_{n-k}^{p^k}\otimes\tau_k+\tau_n\otimes1$$
The conjugation $\iota:\ca_{p*}\to\ca_{p*}$ is determined by the following equalities. 
(See \cite{[17]} for more explicit formula for $\iota(\xi_n)$.)
$$\sum_{k=0}^n\zeta_{n-k}^{2^k}\iota(\zeta_k)=0, \qquad 
\sum_{i=0}^n\xi_{n-k}^{p^k}\iota(\xi_k)=0, \qquad
\iota(\tau_n)=-\sum\limits_{k=0}^n\iota(\xi_{n-k})^{p^k}\tau_k$$

We denote by $\alg_{\bff_p*}$ the category of graded commutative algebras over a prime field 
$\bff_p$ and by $\set$ the category of sets. 
For an object $R_*$ of $\alg_{\bff_p*}$, let $h_{R_*}:\alg_{\bff_p*}\to\set$ be the functor defined 
as follows. 
Put $h_{R_*}(A_*)=\alg_{\bff_p*}(R_*,A_*)=(\text{the set of all homomorphisms from $R_*$ to $A_*$})$
for an object $A_*$ of $\alg_{\bff_p}$ and let $h_{R_*}(f):h_{R_*}(A_*)\to h_{R_*}(B_*)$ map 
$\alpha\in h_{R_*}(A_*)$ to $f\alpha\in h_{R_*}(B_*)$ for a homomorphism $f:A_*\to B_*$ of 
$\alg_{\bff_p*}$. 
If $R_*$ has a structure of Hopf algebra, the coproduct of $R_*$ defines a group structure on 
$h_{R_*}(A_*)$, hence $h_{R_*}$ takes values in the category of groups. 
In this case, $h_{R_*}$ is called the affine group scheme represented by $R_*$ (See Chapter 1 of 
\cite{WCW}.).  
It is shown in the appendix of \cite{ay1} that the affine group scheme represented by $\ca_{p*}$ is 
naturally equivalent to a $\bff_p$-group functor $G_p$ (See Definition \ref{steenrod group} for the 
precise definition.) which assigns each 
graded commutative algebra $A_*$ over $\bff_p$ to a group of isomorphisms of the additive formal 
group law over $A_*[\epsilon]/(\epsilon^2)$ whose coefficient $\alpha_0$ of the leading term 
satisfies $\alpha_0-1\in(\epsilon)$ if $p$ is an odd prime. 

Milnor defined ``The homomorphism $\lambda^*$'' in the section 4 of \cite{[17]} which played a 
crucial role to determine the structure of $\ca_{p*}$. 
This homomorphism $\lambda^*$ also gives a right $\ca_{p*}$-comodule structure on the mod $p$ 
cohomology of a finite complex. 
Hence the mod $p$ cohomology of a finite complex is regarded as a representation of $G_p$ 
(\cite{WCW}, Chapter 3). 
In general, it can be shown that the category of left $\ca_p^*$-modules of finite type is isomorphic 
to the category of right $\ca_{p*}$-comodule of finite type (Theorem 5.2.9 of \cite{rsg} under more 
general settings) and the latter category is isomorphic to the category of representations of $G_p$. 
Thus the theory of modules over the Streenrod algebra is nothing but the representation theory of 
the group scheme $G_p$, which motivates to study the group structure of $G_p$.  

In this paper, we investigate some group theoretical properties of $G_p$. 
Namely, in section 3, we define a descending filtration 
$$G_p\supset G_p^{(0)}\supset G_p^{(0.5)}\supset G_p^{(1)}\supset\cdots\supset 
G_p^{(k)}\supset G_p^{(k+0.5)}\supset G_p^{(k+1)}\supset\cdots$$ 
of normal subgroup schemes of $G_p$ and show that $\Gamma_{k+1}(G_p(A_*))\subset G_p^{(k+0.5)}(A_*)$ 
for a non-negative integer $k$ and a graded commutative algebra $A_*$ over $\bff_p$ 
(\ref{commutator cor}). 
Here $\Gamma_k(G)$ is a subgroup of a group $G$ defined inductively by $\Gamma_0(G)=G$ and 
$\Gamma_{k+1}(G)=[\,\Gamma_k(G),G\,]$. 
Let $\ca_2(n)^*$ be the Hopf subalgebra of $\ca_2^*$ generated by $\sq^{2^i}$ for $i=0,1,\dots,n-1$ 
and $\ca_p(n)^*$ the Hopf subalgebra of $\ca_p^*$ generated by $\beta$ and $\wp^{p^i}$ for 
$i=0,1,\dots,n-1$ if $p$ is an odd prime. 
It follows from Proposition 2 of \cite{[17]} that the dual Hopf algebra $\ca_p(n)_*$ of $\ca_p(n)^*$ 
is a quotient Hopf algebra of $\ca_{p*}$ given by 
$$\ca_p(n)_*=\begin{cases}
\ca_{2*}/\bigl(\zeta_1^{2^n},\dots,\zeta_i^{2^{n-i+1}},\dots,\zeta_n^2,\zeta_{n+1},\zeta_{n+2}, \dots
\bigr)&\text{if $p=2$}\\
\ca_{p*}/\bigl(\xi_1^{p^n},\dots,\xi_i^{p^{n-i+1}},\dots,\xi_n^p,
\tau_{n+1},\xi_{n+1},\tau_{n+2},\xi_{n+2},\dots\bigr)&\text{if $p$ is an odd prime.}\end{cases}$$
We consider a finite subgroup scheme $G_{p,n}$ of $G_p$ represented by $\ca_p(n)_*$ and show that 
the nilpotency class of $G_{p,n}(A_*)$ is at most $n+1$ for any graded commutative $\bff_p$-algebra 
$A_*$. 
In section 4, we give the following diagram in the category of group schemes over $\bff_p$ such that 
$1\to G_p^{ab}\xrightarrow{\kappa_0}G_p\xrightarrow{\rho_0}G_p^{\langle1\rangle}\to1$ and  
$1\to G_p^{[k]}\xrightarrow{\kappa_k}G_p^{\langle k\rangle}\xrightarrow{\rho_k}
G_p^{\langle k+1\rangle}\to1$ are short exact sequences of group schemes and that 
$G_p^{ab}$ and $G_p^{[k]}$ are ``maximal abelian subgroup schemes'' of $G_p$ and 
$G_p^{\langle k\rangle}$ in the sense that $G_p^{ab}$ and $G_p^{[k]}$ are closed subgroup schemes 
defined by minimal Hopf ideals of the Hopf algebras representing $G_p$ and 
$G_p^{\langle k\rangle}$ such that $G_p^{ab}$ and $G_p^{[k]}$ are abelian. 

$$\begin{CD}
G_p^{ab} @. G_p^{[1]} @. G_p^{[2]} @.\cdots @. G_p^{[k]} @. G_p^{[k+1]} @. \cdots \\
@VV{\kappa_0}V @VV{\kappa_1}V @VV{\kappa_2}V @. @VV{\kappa_k}V @VV{\kappa_{k+1}}V \\
G_p @>{\rho_0}>> G_p^{\langle1\rangle} @>{\rho_1}>> G_p^{\langle2\rangle} @>{\rho_2}>> \cdots 
@>{\rho_{k-1}}>> G_p^{\langle k\rangle} @>{\rho_k}>> G_p^{\langle k+1\rangle} @>{\rho_{k+1}}>> \cdots
\end{CD}$$

\section{The Steenrod group}

\begin{Def}\label{steenrod group}
We define $\bff_p$-group functors $G_p$, $G_p^{ev}$ and $G_p^{od}$ as follows. 
Let $A_*$ be a graded commutative algebra over $\bff_p$. 
If $p=2$, we assign degree $-1$ to a variable $X$ and define $G_2(A_*)$ to be the following 
subset of $A_*[[X]]$. 
$$\biggl\{\alpha(X)\in A_*[[X]]\,\biggl|\,\alpha(X)=\sum\limits_{i=0}^{\infty}\alpha_iX^{2^i},\;
\deg\,\alpha_i=2^i-1\;(i\geqq0),\,\alpha_0=1\biggr\}$$
If $p$ is an odd prime, we assign degree $-2$ to a variable $X$ and consider a graded exterior 
algebra $\bff_p[\epsilon]/(\epsilon^2)$ with $\deg\,\epsilon=-1$. 
Define $G_p(A_*)$ to be the following subset of 
$A_*\!\otimes_{\bff_p}\!\bff_p[\epsilon]/(\epsilon^2)[[X]]=A_*[\epsilon]/(\epsilon^2)[[X]]$. 
$$\biggl\{\alpha(X)\in A_*[\epsilon]/(\epsilon^2)[[X]]\,\biggl|\,\alpha(X)=\sum\limits_{i=0}^{\infty}
\alpha_iX^{p^i},\deg\,\alpha_i=2(p^i-1)\;(i\geqq0),\,\alpha_0-1\in(\epsilon)\biggr\}$$
We give a group structure to $G_p(A_*)$ by the composition of formal power series. 
Namely, the product $\alpha(X)\cdot\beta(X)$ of $\alpha(X)$ and $\beta(X)$ is defined by 
$\alpha(X)\cdot\beta(X)=\beta(\alpha(X))$. 
We call $G_p$ the mod $p$ Steenrod group. 
For an odd prime $p$, $G_p^{ev}(A_*)$ is defined to be the following subset of $A_*[[X]]$. 
$$\biggl\{\alpha(X)\in A_*[[X]]\,\biggl|\,\alpha(X)=\sum\limits_{i=0}^{\infty}\alpha_iX^{p^i},
\deg\,\alpha_i=2(p^i-1)\;(i\geqq0),\,\alpha_0=1\biggr\}$$
Since $A_*$ is a subalgebra of $A_*[\epsilon]/(\epsilon^2)$, we regard $G_p^{ev}(A_*)$ as a 
subgroup of $G_p(A_*)$. 
We also define $G_p^{od}(A_*)$ to be the following subset of $A_*[\epsilon]/(\epsilon^2)[[X]]$. 
$$\biggl\{\alpha(X)\in A_*[\epsilon]/(\epsilon^2)[[X]]\,\biggl|\,\alpha(X)=\sum\limits_{i=0}^{\infty}
\alpha_iX^{p^i},\deg\,\alpha_i=2(p^i-1)\;(i\geqq0),\,\alpha_0-1,\alpha_i\in(\epsilon)\;(i\geqq1)
\biggr\}$$
\end{Def}

\begin{Remark}\label{steenrod group rem}
$(1)$ For $\alpha(X)=\sum\limits_{i=0}^{\infty}\alpha_iX^{p^i}, 
\beta(X)=\sum\limits_{i=0}^{\infty}\beta_iX^{p^i}\in G_p(X)$, 
the product $\alpha(X)\cdot\beta(X)$ of $\alpha(X)$ and $\beta(X)$ is given as follows.  
$$\alpha(X)\cdot\beta(X)=
\sum\limits_{i=0}^{\infty}\beta_i\biggl(\sum\limits_{j=0}^{\infty}\alpha_jX^{p^j}\biggr)^{p^i}=
\sum\limits_{i=0}^{\infty}\sum\limits_{j=0}^{\infty}\alpha_j^{p^i}\beta_iX^{p^{i+j}}=
\sum\limits_{i=0}^{\infty}\biggl(\sum\limits_{j=0}^i\alpha_{i-j}^{p^j}\beta_j\biggr)X^{p^i}$$

$(2)$ The quotient map $q_{A_*}:A_*[\epsilon]/(\epsilon^2)\to A_*[\epsilon]/(\epsilon)=A_*$ defines 
a homomorphism 
$$\pi_{A_*}^{ev}:G_p(A_*)\to G_p^{ev}(A_*)$$ 
of groups which is a left inverse of the inclusion map and that $G_p^{od}(A_*)$ is the kernel of 
$\pi_{A_*}^{ev}$. 
Hence $G_p(A_*)$ is a semi-direct product of $G_p^{od}(A_*)$ and $G_p^{ev}(A_*)$. 
\end{Remark}

We denote by $\ca_{p*}^{ev}$ the polynomial part $\bff_p[\xi_1,\xi_2,\dots]$ of $\ca_{p*}$ 
which is a Hopf subalgebra of $\ca_{p*}$. 

\begin{Prop}\label{steenrod group2}(G. Nishida,\;\cite{ay1})
The mod $p$ dual Steenrod algebra $\ca_{p*}$ represents $G_p$. 
If $p$ is an odd prime, $\ca_{p*}^{ev}$ represents $G_p^{ev}$. 
\end{Prop}

\begin{Remark}\label{steenrod group2 rem}
We denote by $\ca_{p*}^{od}$ the quotient of $\ca_{p*}$ by the ideal generated by 
$\xi_1,\xi_2,\dots,\xi_n,\dots$. 
Then, we have $\ca_{p*}^{od}=E(\tau_0,\tau_1,\dots,\tau_n,\dots)$ and each $\tau_i$ is primitive. 
$\ca_{p*}^{od}$ represents $G_p^{od}$ and the quotient map $\ca_{p*}\to\ca_{p*}^{od}$ induces the 
inclusion morphism $G_p^{od}\to G_p$. 
In fact, for a graded commutative $\bff_p$-algebra $A_*$, the natural bijection is given by 
assigning a morphism $f:\ca_{p*}^{od}\to A_*$ to an element 
$(1+f(\tau_0)\epsilon)X+\sum\limits_{i=1}^{\infty}f(\tau_i)\epsilon X^{p^i}$ of $G_p^{od}(A_*)$. 
Since $\alpha(X)^p=X^p$ if $\alpha(X)\in G_p^{od}(A_*)$, $G_p^{od}(A_*)$ is an abelian $p$-subgroup 
of $G_p(A_*)$. 
In fact, $G_p^{od}(A_*)$ is isomorphic to $\prod\limits_{i=0}^{\infty}A_{2p^i-1}$ as an additive 
group. 
\end{Remark}

For a positive integer $n$, an ordered partition of $n$ is a sequence $(\nu(1),\nu(2),\dots,\nu(l))$ 
of positive integers which satisfies $\nu(1)+\nu(2)+\cdots+\nu(l)=n$. 
We denote by $\Part(n)$ the set of all partitions of $n$. 
For a partition $\nu=(\nu(1),\nu(2),\dots,\nu(l))$ of $n$, we put $\ell(\nu)=l$ and 
$\sigma(\nu)(i)=\sum\limits_{s=1}^{i-1}\nu(s)$ ($1\leqq i\leqq l$). 
We call $\ell(\nu)$ the length of $\nu$ and denote by $\Part_l(n)$ the subset of $\Part(n)$ 
consisting of partitions of length $l$. 

\begin{Lemma}\label{partition}
For integers $1\leqq l\leqq k<m$, we define a map $F_{l,k}:\Part_l(k)\to\Part_{l+1}(m)$ by 
$$F_{l,k}((\nu(1),\nu(2),\dots,\nu(l)))=(\nu(1),\nu(2),\dots,\nu(l),m-k).$$ 
Let $F:\bigcup\limits_{k=1}^{m-1}\Part(k)\to\Part(m)$ be the map induced by $F_{l,k}$'s. 
Then, $F$ is an injection whose image is partitions of $m$ of length greater than one. 
\end{Lemma}

\begin{Proof}
Since each $F_{l,k}$ is injective and the images of $F_{l,k}$'s are disjoint each other, $F$ is 
injective. 
For each $\nu=(\nu(1),\nu(2),\dots,\nu(l),\nu(l+1))\in\Part_{l+1}(m)$, $F_{l,\sigma(\nu)(l+1)}$ 
maps $(\nu(1),\nu(2),\dots,\nu(l))\in\Part_l(\sigma(\nu)(l+1))$ to $\nu$.  
\end{Proof}

\begin{Prop}\label{anti-auto}(\cite{[17]})
Let $A_*$ be a graded commutative algebra over $\bff_p$ and $c$ a fixed integer which is even if 
$p$ is odd. 
Suppose that sequences of elements $(\alpha_i)_{i\geqq0}$ and $(\beta_i)_{i\geqq0}$ of $A_*$ satisfy 
the following conditions. 
\begin{description}
\item{$(i)$} $\alpha_0=1$ if $p=2$, $(\alpha_0-1)^2=0$ if $p$ is odd. 
\item{$(ii)$} $\deg\,\alpha_0=\deg\,\beta_0=0$, $\deg\,\alpha_i=\deg\,\beta_i=c(1+p+\cdots+p^{i-1})$. 
\item{$(iii)$} $\alpha_0\beta_0=1$ and $\sum\limits_{k=0}^i\alpha_{i-k}^{p^k}\beta_k=0$ for any 
positive integer $i$. 
\end{description}
Then, $\beta_0=\alpha_0^{-1}=2-\alpha_0$ and the following equality holds for each positive integer 
$n$. 
$$\beta_n=\beta_0\sum\limits_{\nu\in\Part(n)}(-1)^{\ell(\nu)}\prod\limits_{j=1}^{\ell(\nu)}
\alpha_{\nu(j)}^{p^{\sigma(\nu)(j)}}$$
\end{Prop}

\begin{Proof}
We have $\alpha_0(2-\alpha_0)=1$ by $(i)$. 
Thus it follows from $\alpha_0\beta_0=1$ that $\beta_0=\alpha_0^{-1}=2-\alpha_0$. 
We put $\tilde\beta_i=\beta_i\alpha_0$. 
Then $\tilde\beta_0=1$ and $\sum\limits_{k=0}^i\alpha_{i-k}^{p^k}\tilde\beta_k=0$ for any positive 
integer $i$ by $(iii)$. 
It suffices to show the following.  
$$\tilde\beta_n=\sum\limits_{\nu\in\Part(n)}(-1)^{\ell(\nu)}\prod\limits_{j=1}^{\ell(\nu)}
\alpha_{\nu(j)}^{p^{\sigma(\nu)(j)}}$$
Since $\alpha_0=1+(\alpha_0-1)$ and $(\alpha_0-1)^2=0$, we have $\alpha_0^p=1$. 
Then, $\alpha_1\tilde\beta_0+\alpha_0^p\tilde\beta_1=0$ implies $\tilde\beta_1=-\alpha_1$ and the 
assertion holds for $n=1$. 
Suppose that assertion holds for $1\leqq n\leqq m-1$. 
We consider the map $F:\bigcup\limits_{k=1}^{m-1}\Part(k)\to\Part(m)$ in (\ref{partition}). 
For $\nu\in\Part(k)$, we have $\sigma(F(\nu))(j)=\sigma(\nu)(j)$ if $1\leqq j\leqq\ell(\nu)$ and 
$\sigma(F(\nu))(\ell(\nu)+1)=k$. 
Hence it follows from the inductive hypothesis and (\ref{partition}) that 
\begin{align*}
\tilde\beta_m&=-\sum\limits_{k=0}^{m-1}\alpha_{m-k}^{p^k}\tilde\beta_k=
-\alpha_m\tilde\beta_0-\sum\limits_{k=1}^{m-1}\alpha_{m-k}^{p^k}\sum\limits_{\nu\in\Part(k)}
(-1)^{\ell(\nu)}\prod\limits_{j=1}^{\ell(\nu)}\alpha_{\nu(j)}^{p^{\sigma(\nu)(j)}}\\
&=-\alpha_m+\sum\limits_{k=1}^{m-1}\sum\limits_{\nu\in\Part(k)}(-1)^{\ell(F(\nu))}
\alpha_{F(\nu)(\ell(\nu)+1)}^{p^{\sigma(F(\nu))(\ell(\nu)+1)}}
\prod\limits_{j=1}^{\ell(\nu)}\alpha_{F(\nu)(j)}^{p^{\sigma(F(\nu))(j)}}\\
&=-\alpha_m+\sum\limits_{k=1}^{m-1}\sum\limits_{\nu\in\Part(k)}(-1)^{\ell(F(\nu))}
\prod\limits_{j=1}^{\ell(\nu)+1}\alpha_{F(\nu)(j)}^{p^{\sigma(F(\nu))(j)}}=
\sum\limits_{\nu\in\Part(m)}(-1)^{\ell(\nu)}\prod\limits_{j=1}^{\ell(\nu)}
\alpha_{\nu(j)}^{p^{\sigma(\nu)(j)}}.
\end{align*}
\end{Proof}

The next result is a direct consequence of (\ref{anti-auto}) and the above equality. 

\begin{Prop}\label{inverse}
The inverse of $\alpha(X)=\sum\limits_{i=0}^{\infty}\alpha_iX^{p^i}\in G_p(A_*)$ is given as follows.  
$$\alpha(X)^{-1}=\alpha_0^{-1}X+\sum\limits_{i=1}^{\infty}\alpha_0^{-1}\!\left(
\sum\limits_{\nu\in\Part(i)}(-1)^{\ell(\nu)}\prod\limits_{j=1}^{\ell(\nu)}
\alpha_{\nu(j)}^{p^{\sigma(\nu)(j)}}\right)\!X^{p^i}$$
\end{Prop}

\begin{Remark}\label{inverse rem}
Suppose that $p$ is an odd prime and that $\alpha_i=\alpha_{0i}+\alpha_{1i}\epsilon$ for 
$\alpha_{0i}\in A_{2(p^i-1)}$, $\alpha_{1i}\in A_{2p^i-1}$. 
For $\nu\in\Part(i)$, we have 
$$\prod\limits_{j=1}^{\ell(\nu)}\alpha_{\nu(j)}^{p^{\sigma(\nu)(j)}}=
\bigl(\alpha_{0\nu(1)}+\alpha_{1\nu(1)}\epsilon\bigr)
\prod\limits_{j=2}^{\ell(\nu)}\alpha_{0\nu(j)}^{p^{\sigma(\nu)(j)}}=
\prod\limits_{j=1}^{\ell(\nu)}\alpha_{0\nu(j)}^{p^{\sigma(\nu)(j)}}+\left(\alpha_{1\nu(1)}
\prod\limits_{j=2}^{\ell(\nu)}\alpha_{0\nu(j)}^{p^{\sigma(\nu)(j)}}\right)\!\epsilon.$$
Hence we have the following formula. 
$$\alpha(X)^{-1}=(1-\alpha_{10}\epsilon)X+\sum\limits_{i=1}^{\infty}(1-\alpha_{10}\epsilon)\!\left(
\sum\limits_{\nu\in\Part(i)}(-1)^{\ell(\nu)}\!\left(\prod\limits_{j=1}^{\ell(\nu)}
\alpha_{0\nu(j)}^{p^{\sigma(\nu)(j)}}+\left(\alpha_{1\nu(1)}\prod\limits_{j=2}^{\ell(\nu)}
\alpha_{0\nu(j)}^{p^{\sigma(\nu)(j)}}\right)\!\epsilon\right)\right)\!X^{p^i}$$
\end{Remark}

\section{Lower central series}

For a non-negative integer $k$, we define ``quotient groups" $G_p^k$ of $G_p$ as follows. 
If $p=2$, we assign degree $-1$ to a variable $X$ and define $G_2^k(A_*)$ to be the following 
subset of $A_*[X]/(X^{2^{k+1}})$. 
$$\biggl\{\alpha(X)\in A_*[X]/(X^{2^{k+1}})\,\biggl|\,\alpha(X)=\sum\limits_{i=0}^k\alpha_iX^{2^i},\;
\deg\,\alpha_i=2^i-1\,(0\leqq i\leqq k),\,\alpha_0=1\biggr\}$$
If $p$ is an odd prime, we assign degree $-2$ to $X$ and consider a graded exterior 
algebra $\bff_p[\epsilon]/(\epsilon^2)$ with $\deg\,\epsilon=-1$. 
Let $G_p^k(A_*)$ to be a subset of 
$A_*\!\otimes_{\bff_p}\!\bff_p[\epsilon]/(\epsilon^2)[X]/(X^{p^{k+1}})=
A_*[\epsilon]/(\epsilon^2)[X]/(X^{p^{k+1}})$ given by 
$$\biggl\{\alpha(X)\in A_*[\epsilon]/(\epsilon^2)[X]/(X^{p^{k+1}})\,\biggl|\,
\alpha(X)=\sum\limits_{i=0}^k\alpha_iX^{p^i},\deg\,\alpha_i=2(p^i-1)\,(0\leqq i\leqq k),\,
\alpha_0-1\in(\epsilon)\biggr\}.$$
We give a group structure to $G_p^k(A_*)$ by the composition of stunted polynomials. 
Namely, 
$$\alpha(X)\cdot\beta(X)=\beta(\alpha(X))=
\sum\limits_{i=0}^k\beta_i\biggl(\sum\limits_{j=0}^k\alpha_jX^{p^j}\biggr)^{p^i}=
\sum\limits_{i=0}^k\sum\limits_{j=0}^k\alpha_j^{p^i}\beta_iX^{p^{i+j}}=
\sum\limits_{i=0}^k\biggl(\sum\limits_{l=0}^i\alpha_{i-l}^{p^l}\beta_l\biggr)X^{p^i}.$$

Define maps $\pi^k_{A_*}:G_p(A_*)\to G_p^k(A_*)$ ($k=0,1,2,\dots$) to be the restrictions of the 
quotient maps $A_*[[X]]\to A_*[X]/(X^{2^{k+1}})$ if $p=2$ and $A_*[\epsilon]/(\epsilon^2)[[X]]\to 
(A_*[\epsilon]/(\epsilon^2)[X])/(X^{p^{k+1}})$ if $p$ is an odd prime. 
It is clear that $\pi^k_{A_*}$ is a homomorphism of groups and natural in $A_*$. 
We denote by $G_p^{(k)}(A_*)$ the kernel of $\pi^k_{A_*}$, that is, 
$$G_p^{(k)}(A_*)=\biggl\{\alpha(X)\in G_p(A_*)\,\biggl|\,
\alpha(X)=X+\sum\limits_{i=k+1}^{\infty}\alpha_iX^{p^i}\biggr\}.$$

We regard $A_*$ as a subalgebra of $A_*[\epsilon]/(\epsilon^2)$ and define a subset 
$G_p^{k+0.5}(A_*)$ of $G_p^{k+1}(A_*)$ by   
$$G_p^{k+0.5}(A_*)=\biggl\{\alpha(X)\in G_p^{k+1}(A_*)\,\biggl|\,
\alpha(X)=\sum\limits_{i=0}^{k+1}\alpha_iX^{p^i},\,\alpha_{k+1}\in A_{2p^{k+1}-2}\biggr\}.$$
Define a map $q_{A_*}^k:G_p^{k+1}(A_*)\to G_p^{k+0.5}(A_*)$ by 
$$q_{A_*}^k(\alpha(X))=\sum\limits_{i=0}^k\alpha_iX^{p^i}+q_{A_*}(\alpha_{k+1})X^{p^{k+1}}$$ 
if $\alpha(X)=\sum\limits_{i=0}^{k+1}\alpha_iX^{p^i}$. 
Here $q_{A_*}:A_*[\epsilon]/(\epsilon^2)\to A_*[\epsilon]/(\epsilon)=A_*$ denotes the quotient map 
again.  
For $\alpha(X),\beta(X)\in G_p^{k+0.5}(A_*)$, we set 
$\alpha(X){*}\beta(X)=q_{A_*}^{k+1}(\alpha(X){\cdot}\beta(X))$. 
Then, the correspondence $(\alpha(X),\beta(X))\mapsto\alpha(X){*}\beta(X)$ defines a group structure 
on $G_p^{k+0.5}(A_*)$. 
In fact, the inverse of $\alpha(X)$ is $q_{A_*}^{k+1}(\alpha(X)^{-1})$. 
We note that $q_{A_*}^k$ is a homomorphism of groups. 
Let us denote by $G_p^{(k+0.5)}(A_*)$ the kernel of a composition 
$G_p(A_*)\xrightarrow{\pi^{k+1}_{A_*}}G_p^{k+1}(A_*)\xrightarrow{q_{A_*}^k}G_p^{k+0.5}(A_*)$.  
Then, we have 
$$G_p^{(k+0.5)}(A_*)=\biggl\{\alpha(X)\in G_p(A_*)\,\bigg|\,
\alpha(X)=X+\sum\limits_{i=k+1}^{\infty}\alpha_iX^{p^i},\alpha_{k+1}\in(\epsilon)\biggr\}$$ 
Here we put $\epsilon=0$ if $p=2$. 
Then, $G_2(A_*)=G_2^{(0)}(A_*)$, $G_2^{(k+0.5)}(A_*)=G_2^{(k+1)}(A_*)$ and we have a 
decreasing filtration of $G_p(A_*)$. 
$$G_p(A_*)\supset G_p^{(0)}(A_*)\supset G_p^{(0.5)}(A_*)\supset G_p^{(1)}(A_*)\supset\cdots\supset 
G_p^{(k)}(A_*)\supset G_p^{(k+0.5)}(A_*)\supset G_p^{(k+1)}(A_*)\supset\cdots$$

\begin{Lemma}\label{filtration of rsg}
Suppose that elements $\alpha(X)=\sum\limits_{i=0}^{\infty}\alpha_iX^{p^i}$ and 
$\beta(X)=\sum\limits_{i=0}^{\infty}\beta_iX^{p^i}$ of $G_p(A_*)$ satisfy $\alpha_i=0$ for 
$i=1,2,\dots,k$ and $\beta_i=0$ for $i=1,2,\dots,l$, respectively. 
We put $\beta(X)^{-1}=\sum\limits_{i=0}^{\infty}\bar\beta_iX^{p^i}$. 

$(1)$ If $k=l=0$, the following equality holds. 
\begin{align*}
[\alpha(X),\beta(X)]&=
X+(\alpha_1(\beta_0-1)+(1-\alpha_0)\beta_1)X^p\\
&{\phantom =}\;+((\alpha_2-\alpha_1^{p+1})(\beta_0-1)+
(1-\alpha_0)(\beta_2-\beta_1^{p+1})+\alpha_0\alpha_1^p\beta_1-\alpha_1\beta_0\beta_1^p)X^{p^2}\\
&{\phantom =}\;+(\text{higher terms})
\end{align*}

$(2)$ If $k\geqq1$ and $l=0$, the following equality holds. 
\begin{align*}
[\alpha(X),\beta(X)]&=
X+(\alpha_{k+1}(\beta_0-1)-(1-\alpha_0)\beta_0\bar\beta_{k+1})X^{p^{k+1}}\\
&{\phantom =}\;+\bigl(\alpha_{k+2}(\beta_0-1)-(1-\alpha_0)\beta_0\bar\beta_{k+2}+
\alpha_0\alpha_{k+1}^p\beta_1-\alpha_{k+1}\beta_0\beta_1^{p^{k+1}}\bigr)X^{p^{k+2}}\\
&{\phantom =}\;+(\text{higher terms})
\end{align*}

$(3)$ If $k\geqq l\geqq1$, the following equality holds. 
\begin{align*}
[\alpha(X),\beta(X)]&=
X+(\alpha_{k+1}(\beta_0-1)-(1-\alpha_0)\beta_0\bar\beta_{k+1})X^{p^{k+1}}+
(\alpha_{k+2}(\beta_0-1)-(1-\alpha_0)\beta_0\bar\beta_{k+2})X^{p^{k+2}}\\
&{\phantom =}\;+(\text{higher terms})
\end{align*}
\end{Lemma}

\begin{Proof}
Since $\alpha_0^p=\beta_0^p=1$ and $\alpha(X)\cdot\beta(X)=
\sum\limits_{i=0}^{\infty}\biggl(\sum\limits_{j=0}^i\alpha_{i-j}^{p^j}\beta_j\biggr)X^{p^i}$, we 
have the following equality. 
$$\alpha(X)\cdot\beta(X)=
\alpha_0\beta_0X+\sum\limits_{i=1}^{k+1}(\alpha_i\beta_0+\beta_i)X^{p^i}\!+
(\alpha_{k+2}\beta_0+\alpha_{k+1}^p\beta_1+\beta_{k+2})X^{p^{k+2}}\!+(\text{higher terms})$$
Hence if we put $\alpha(X)\cdot\beta(X)=\sum\limits_{i=0}^{\infty}\gamma_iX^{p^i}$, $\gamma_i$'s  
are given by $\gamma_0=\alpha_0\beta_0$, $\gamma_i=\beta_i$ for $1\leqq i\leqq k$, 
$\gamma_{k+1}=\alpha_{k+1}\beta_0+\beta_{k+1}$ and 
$\gamma_{k+2}=\alpha_{k+2}\beta_0+\alpha_{k+1}^p\beta_1+\beta_{k+2}$. 

Put $\alpha(X)^{-1}=\sum\limits_{i=0}^{\infty}\bar\alpha_iX^{p^i}$.  
If $1\leqq i\leqq k$ and $\nu\in\Part(i)$, then $\nu(j)\leqq k$ for any $1\leqq j\leqq\ell(\nu)$. 
Since $\alpha_i=0$ for $1\leqq i\leqq k$, we have $\bar\alpha_i=0$ for $1\leqq i\leqq k$ by 
(\ref{inverse}). 
If $\nu\in\Part(k+1)$ satisfies $\nu(j)\geqq k+1$ for $1\leqq j\leqq\ell(\nu)$, then 
$\nu=(k+1)$, hence (\ref{inverse}) implies $\bar\alpha_{k+1}=-\alpha_0^{-1}\alpha_{k+1}$. 
If $\nu\in\Part(k+2)$ satisfies $\nu(j)\geqq k+1$ for $1\leqq j\leqq\ell(\nu)$, then 
$\nu=(k+2)$ or ``$\nu=(1,1)$ and $k=0$''. 
It follows from (\ref{inverse}) that $\bar\alpha_{k+2}=\begin{cases}\alpha_0^{-1}
(\alpha_1^{p+1}-\alpha_2)&k=0\\\hfill-\alpha_0^{-1}\alpha_{k+2}\hfill&k\geqq1\end{cases}$ holds.
Similarly, we have $\bar\beta_1=-\beta_0^{-1}\beta_1$ and 
$\bar\beta_2=\beta_0^{-1}(\beta_1^{p+1}-\beta_2)$ if $l=0$. 
Hence, if we put $\alpha(X)^{-1}\!\cdot\beta(X)^{-1}=\sum\limits_{i=0}^{\infty}\bar\gamma_iX^{p^i}$, 
then $\bar\gamma_i$'s are given by $\bar\gamma_0=\alpha_0^{-1}\beta_0^{-1}$, 
$\bar\gamma_i=\bar\beta_i$ for $1\leqq i\leqq k$ and  
\begin{align*}
\bar\gamma_{k+1}&=\bar\alpha_{k+1}\beta_0^{-1}+\bar\beta_{k+1}=
\bar\beta_{k+1}-\alpha_0^{-1}\alpha_{k+1}\beta_0^{-1}
\\
\bar\gamma_{k+2}&=\bar\alpha_{k+2}\beta_0^{-1}+\bar\alpha_{k+1}^p\bar\beta_1+\bar\beta_{k+2}=
\begin{cases}
\beta_0^{-1}(\alpha_0^{-1}(\alpha_1^{p+1}-\alpha_2)+\alpha_1^p\beta_1+\beta_1^{p+1}-\beta_2)&k=l=0
\\[2mm]
-\alpha_0^{-1}\alpha_{k+2}\beta_0^{-1}+\alpha_{k+1}^p\beta_0^{-1}\beta_1+\bar\beta_{k+2}
&k\geqq1,\,l=0
\\[2mm]
-\alpha_0^{-1}\alpha_{k+2}\beta_0^{-1}+\bar\beta_{k+2}&k\geqq l\geqq1.\end{cases}
\end{align*} 
It follows that $\bar\gamma_0\gamma_0=\alpha_0^{-1}\beta_0^{-1}\alpha_0\beta_0=1$ and  
$\sum\limits_{j=0}^i\bar\gamma_{i-j}^{p^j}\gamma_j=
\sum\limits_{j=0}^i\bar\beta_{i-j}^{p^j}\beta_j=0$ if $1\leqq i\leqq k$. 
We also have 
\begin{align*}
\sum\limits_{j=0}^{k+1}\bar\gamma_{k+1-j}^{p^j}\gamma_j&=
\bar\gamma_{k+1}\gamma_0+\sum\limits_{j=1}^k\bar\gamma_{k+1-j}^{p^j}\gamma_j+
\bar\gamma_0^{p^{k+1}}\gamma_{k+1}=
\alpha_{k+1}(\beta_0-1)-(1-\alpha_0)\beta_0\bar\beta_{k+1}+
\sum\limits_{j=0}^{k+1}\bar\beta_{k+1-j}^{p^j}\beta_j\\
&=\alpha_{k+1}(\beta_0-1)-(1-\alpha_0)\beta_0\bar\beta_{k+1}
\end{align*}
If $k=l=0$, we have 
$$\sum\limits_{j=0}^2\bar\gamma_{2-j}^{p^j}\gamma_j
=\bar\gamma_2\gamma_0+\bar\gamma_1^p\gamma_1+\bar\gamma_0^{p^2}\gamma_2
=(\alpha_2-\alpha_1^{p+1})(\beta_0-1)+(1-\alpha_0)(\beta_2-\beta_1^{p+1})
+\alpha_0\alpha_1^p\beta_1-\alpha_1\beta_0\beta_1^p.$$ 
If $k\geqq1$ and $l=0$, the following equalities hold. 
\begin{align*}
\sum\limits_{j=0}^{k+2}\bar\gamma_{k+2-j}^{p^j}\gamma_j
&=\bar\gamma_{k+2}\gamma_0+\bar\gamma_{k+1}^p\gamma_1+
\sum\limits_{j=2}^k\bar\gamma_{k+2-j}^{p^j}\gamma_j+
\bar\gamma_1^{p^{k+1}}\gamma_{k+1}+\bar\gamma_0^{p^{k+2}}\gamma_{k+2}
\\
&=\alpha_{k+2}(\beta_0-1)-(1-\alpha_0)\beta_0\bar\beta_{k+2}+\alpha_0\alpha_{k+1}^p\beta_1
-\alpha_{k+1}\beta_0\beta_1^{p^{k+1}}+\sum\limits_{j=0}^{k+2}\bar\beta_{k+2-j}^{p^j}\beta_j
\\
&=\alpha_{k+2}(\beta_0-1)-(1-\alpha_0)\beta_0\bar\beta_{k+2}+\alpha_0\alpha_{k+1}^p\beta_1
-\alpha_{k+1}\beta_0\beta_1^{p^{k+1}}
\end{align*}
If $k\geqq l\geqq1$, the following equalities hold. 
\begin{align*}
\sum\limits_{j=0}^{k+2}\bar\gamma_{k+2-j}^{p^j}\gamma_j
&=\bar\gamma_{k+2}\gamma_0+\bar\gamma_{k+1}^p\gamma_1+
\sum\limits_{j=2}^k\bar\gamma_{k+2-j}^{p^j}\gamma_j+
\bar\gamma_1^{p^{k+1}}\gamma_{k+1}+\bar\gamma_0^{p^{k+2}}\gamma_{k+2}
\\
&=\alpha_{k+2}(\beta_0-1)-(1-\alpha_0)\beta_0\bar\beta_{k+2}
+\sum\limits_{j=0}^{k+2}\bar\beta_{k+2-j}^{p^j}\beta_j
=\alpha_{k+2}(\beta_0-1)-(1-\alpha_0)\beta_0\bar\beta_{k+2}
\end{align*}
\end{Proof}

\begin{Prop}\label{commutator}
The following relations hold. 
\begin{align*}
[\,G_p(A_*),G_p(A_*)\,]&\subset G_p^{(0.5)}(A_*)\\
[\,G_p^{(0.5)}(A_*),G_p^{(0.5)}(A_*)\,]&\subset G_p^{(2)}(A_*)\\
[\,G_p^{(k)}(A_*),G_p^{(k)}(A_*)\,]&\subset G_p^{(k+2)}(A_*)\;\;\text{if $k$ is a positive integer.}
\\
[\,G_p^{(k)}(A_*),G_p(A_*)\,]&\subset G_p^{(k+0.5)}(A_*)\;\;\text{if $k$ is a non-negative integer.}
\\
[\,G_p^{(k+0.5)}(A_*),G_p(A_*)\,]&\subset G_p^{(k+1.5)}(A_*)\;\;
\text{if $k$ is a non-negative integer.}
\end{align*}
\end{Prop}

\begin{Proof}
The first and the second relations are direct consequences of (1) of (\ref{filtration of rsg}). 
The third relation follows from (3) of (\ref{filtration of rsg}). 
For $\alpha(X)=X+\sum\limits_{i=k+1}^{\infty}\alpha_iX^{p^i}\in G_p^{(k)}(A_*)$ and 
$\beta(X)=\sum\limits_{i=0}^{\infty}\beta_iX^{p^i}\in G_p(A_*)$, since $\beta_0-1\in(\epsilon)$, 
the fourth relation follows from (2) of (\ref{filtration of rsg}). 
If $\alpha(X)\in G_p^{(k+0.5)}(A_*)$, then $\alpha_{k+1}\in(\epsilon)$ which implies that 
$\alpha_{k+1}(\beta_0-1)=0$ and 
$\alpha_{k+2}(\beta_0-1)+\alpha_{k+1}^p\beta_1-\beta_0\alpha_{k+1}\beta_1^{p^{k+1}}\in(\epsilon)$. 
Hence the fifth relation also follows from (2) of (\ref{filtration of rsg}). 
\end{Proof}

For a group $G$ and a non-negative integer $k$, we define subgroups $D_k(G)$ and $\Gamma_k(G)$ of 
$G$ inductively by $D_0(G)=\Gamma_0(G)=G$ and $D_{k+1}(G)=[\,D_k(G),D_k(G)\,]$, 
$\Gamma_{k+1}(G)=[\,\Gamma_k(G),G\,]$. 
The following result is a direct consequence of (\ref{commutator})

\begin{Cor}\label{commutator cor}
We have $D_1(G_p(A_*))=\Gamma_1(G_p(A_*))\subset G_p^{(0.5)}(A_*)$. 
For a positive integer $k$, the following relations hold. 
$$D_{k+1}(G_p(A_*))\subset G_p^{(2k)}(A_*),\qquad\Gamma_{k+1}(G_p(A_*))\subset G_p^{(k+0.5)}(A_*)$$ 
\end{Cor}

\begin{Remark}\label{commutator rem}
If $H$ is a subgroup of $G_p(A_*)$, we have the following relations. 
$$D_1(H)=\Gamma_1(H)\subset G_p^{(0.5)}(A_*)\cap H, \quad D_{k+1}(H)\subset G_p^{(2k)}(A_*)\cap H, 
\quad\Gamma_{k+1}(H)\subset G_p^{(k+0.5)}(A_*)\cap H$$
Since $G_p^{(k+0.5)}(A_*)\cap G_p^{ev}(A_*)\subset G_p^{(k+1)}(A_*)$, we have 
$\Gamma_{k+1}(H)\subset G_p^{(k+1)}(A_*)\cap H$ if $H$ is a subgroup of $G_p^{ev}(A_*)$. 
\end{Remark}

We define another filtration of $G_p$ as follows. 
For a non-negative integer $n$ and a graded commutative algebra $A_*$ over $\bff_p$, let 
$G_{p,n}(A_*)$ be a subset of $G_p(A_*)$ consisting of elements 
$\alpha(X)=\sum\limits_{i=0}^{\infty}\alpha_iX^{p^i}$ which satisfy 
$\alpha_i^{p^{n-i+1}}=0$ for $i=1,2,\dots,n$ and $\alpha_i=0$ for $i\geqq n+1$. 

\begin{Prop}\label{steenrod subgroup}
$G_{p,n}(A_*)$ is a subgroup of $G_p(A_*)$. 
\end{Prop}

\begin{Proof}
Suppose $\alpha(X)=\sum\limits_{i=0}^{\infty}\alpha_iX^{p^i},
\beta(X)=\sum\limits_{i=0}^{\infty}\beta_iX^{p^i}\in G_{p,n}(A_*)$. 
Put $\alpha(X)\cdot\beta(X)=\sum\limits_{i=0}^{\infty}\gamma_iX^{p^i}$, then we have 
$\gamma_i=\sum\limits_{j=0}^i\alpha_{i-j}^{p^j}\beta_j$. 
Since $\alpha_k^{p^{n-k+1}}=\beta_k^{p^{n-k+1}}=0$ for $k=1,2,\dots,n$, it follows 
$\gamma_i^{p^{n-i+1}}=\sum\limits_{j=0}^i\alpha_{i-j}^{p^{n-(i-j)+1}}\beta_j^{p^{n-i+1}}=0$ if 
$1\leqq i\leqq n$. 
Assume that $i\geqq n+1$. 
Since $\alpha_{i-j}^{p^j}=\alpha_{i-j}^{p^{n-(i-j)+1+(i-n-1)}}=0$ for $i-n\leqq j\leqq n$, we have 
$\gamma_i=\sum\limits_{j=i-n}^n\alpha_{i-j}^{p^j}\beta_j=0$. 
Thus $\alpha(X)\cdot\beta(X)\in G_{p,n}(A_*)$. 

We put $\alpha(X)^{-1}=\alpha_0^{-1}X+\sum\limits_{i=1}^{\infty}\delta_iX^{p^i}$. 
Then,  $\delta_i=\alpha_0^{-1}\!\left(\sum\limits_{\nu\in\Part(i)}(-1)^{\ell(\nu)}
\prod\limits_{j=1}^{\ell(\nu)}\alpha_{\nu(j)}^{p^{\sigma(\nu)(j)}}\right)$ by (\ref{inverse}). 
Suppose that $\prod\limits_{j=1}^{\ell(\nu)}\alpha_{\nu(j)}^{p^{\sigma(\nu)(j)+n-i+1}}\ne0$ for 
some $1\leqq i\leqq n$ and $\nu\in\Part(i)$. 
Then, $\sigma(\nu)(j)+n-i+1\leqq n-\nu(j)$ for $1\leqq j\leqq\ell(\nu)$, which implies a 
contradiction $n+1\leqq n$ if $j=\ell(\nu)$.  
Hence we have $\delta_i^{n-i+1}=0$ for $1\leqq i\leqq n$. 
Suppose that $\prod\limits_{j=1}^{\ell(\nu)}\alpha_{\nu(j)}^{p^{\sigma(\nu)(j)}}\ne0$ for 
some $i\geqq n+1$ and $\nu\in\Part(i)$. 
Then, $\nu(j)\leqq n$ and $\sigma(\nu)(j)\leqq n-\nu(j)$ for $1\leqq j\leqq\ell(\nu)$. 
The latter inequality implies $i\leqq n$ if $j=\ell(\nu)$, which contradicts the assumption.  
Hence we have $\delta_i=0$ for $i\geqq n+1$ and $\alpha(X)^{-1}\in G_{p,n}(A_*)$. 
\end{Proof}

Thus we have the following increasing filtration of subgroups of $G_p(A_*)$. 
$$G_{p,0}(A_*)\subset G_{p,1}(A_*)\subset G_{p,2}(A_*)\subset \cdots\subset G_{p,n}(A_*)
\subset G_{p,n+1}(A_*)\subset\cdots\subset G_p(A_*)$$
If $p$ is an odd prime, we define $G_{p,n}^{ev}(A_*)$ by $G_{p,n}^{ev}(A_*)=G_{p,n}(A_*)\cap 
G_p^{ev}(A_*)$, we have the following increasing filtration of subgroups of $G_p^{ev}(A_*)$. 
$$G_{p,0}^{ev}(A_*)\subset G_{p,1}^{ev}(A_*)\subset G_{p,2}^{ev}(A_*)\subset \cdots\subset 
G_{p,n}^{ev}(A_*)\subset G_{p,n+1}^{ev}(A_*)\subset\cdots\subset G_p^{ev}(A_*)$$

We note that $G_{2,0}(A_*)$ and $G_{p,0}^{ev}(A_*)$ are the trivial groups and that $G_{p,0}(A_*)$ 
is isomorphic to the additive group $A_1$. 
Since $G_p^{(n)}(A_*)\cap G_{p,n}(A_*)$ is the trivial group, (\ref{commutator rem}) implies 
the following fact. 

\begin{Thm}\label{nilpotency}
We have the following lower central series. 
$$G_{p,n}(A_*)
\supset\Gamma_1(G_{p,n}(A_*))\supset\cdots\supset\Gamma_i(G_{p,n}(A_*))\supset
\Gamma_{i+1}(G_{p,n}(A_*))\supset\cdots\supset\Gamma_{n+1}(G_{p,n}(A_*))=\{X\}$$
$$G_{p,n}^{ev}(A_*)
\supset\Gamma_1(G_{p,n}^{ev}(A_*))\supset\cdots\supset\Gamma_i(G_{p,n}^{ev}(A_*))\supset
\Gamma_{i+1}^{ev}(G_{p,n}(A_*))\supset\cdots\supset\Gamma_n(G_{p,n}^{ev}(A_*))=\{X\}$$
\end{Thm}

Let $I_{2,n}$ be an ideal of $\ca_{2*}$ generated by $\zeta_1^{2^n},\zeta_2^{2^{n-1}},\dots,
\zeta_n^2$ and $\zeta_i$ for $i\geqq n+1$. 
For an odd prime $p$, let $I_{p,n}$ be an ideal of $\ca_{p*}$ generated by 
$\xi_1^{p^n},\xi_2^{p^{n-1}},\dots,\xi_n^p$ and $\tau_i$, $\xi_i$ for $i\geqq n+1$ and 
$I_{p,n}^{ev}$ an ideal of $\ca_{p*}^{ev}$ generated by 
$\xi_1^{p^n},\xi_2^{p^{n-1}},\dots,\xi_n^p$ and $\xi_i$ for $i\geqq n+1$. 
We put 
\begin{align*}
\ca_2(n)_*&=\ca_{2*}/I_{2,n}=\bff_2[\zeta_1,\zeta_2,\dots,\zeta_n]/
(\zeta_1^{2^n},\zeta_2^{2^{n-1}},\dots,\zeta_n^2)\\
\ca_p(n)_*&=\ca_{p*}/I_{p,n}=E(\tau_0,\tau_1,\dots,\tau_n)\otimes_{\bff_p}\!
\bff_p[\xi_1,\xi_2,\dots,\xi_n]/(\xi_1^{p^n},\xi_2^{p^{n-1}},\dots,\xi_n^p)\\
\ca_p^{ev}(n)_*&=\ca_{p*}^{ev}/I_{p,n}^{ev}=
\bff_p[\xi_1,\xi_2,\dots,\xi_n]/(\xi_1^{p^n},\xi_2^{p^{n-1}},\dots,\xi_n^p).
\end{align*}

We have the following fact from (\ref{steenrod group2}). 

\begin{Prop}\label{steenrod subgroup2}
$G_{p,n}$ is represented by $\ca_p(n)_*$ and $G_{p,n}^{ev}$ is represented by $\ca_p^{ev}(n)_*$. 
\end{Prop}

\section{A maximal abelian subgroup of the Steenrod group}

We define $\bff_p$-group functors $G_p^{\langle k\rangle}$ for $k=1,2,\dots$ as follows. 

Let $A_*$ be a graded commutative $\bff_p$-algebra. 
If $p=2$, we assign degree $-2^k$ to a variable $X$ and define $G_2^{\langle k\rangle}(A_*)$ by 
$$G_2^{\langle k\rangle}(A_*)=\Biggl\{\sum\limits_{i=0}^{\infty}\alpha_iX^{2^i}\in A_*[[X]]\,
\Biggl|\,\deg\,\alpha_i=2^{i+k}-2^k\;(i\geqq0),\;\alpha_0=1\Biggr\}.$$
If $p$ is an odd prime, we assign degree $-2p^k$ to a variable $X$, degree $-1$ to $\epsilon$ and 
define $G_p^{\langle k\rangle}(A_*)$ by
\begin{align*}
G_p^{\langle1\rangle}(A_*)&=\Biggl\{\sum\limits_{i=0}^{\infty}\alpha_iX^{p^i}\!\in 
A_*[\epsilon]/(\epsilon^2)[[X]]\,\Biggl|\,\deg\,\alpha_i\!=2(p^{i+1}\!-p)\;(i\geqq0),\,
\alpha_0-1\in(\epsilon),\,\alpha_i\in A_*\;(i\geqq1)\Biggr\}, 
\\
G_p^{\langle k\rangle}(A_*)&=\Biggl\{\sum\limits_{i=0}^{\infty}\alpha_iX^{p^i}\!\in 
A_*[[X]]\,\Biggl|\,\deg\,\alpha_i\!=2(p^{i+k}\!-p^k)\;(i\geqq0),\,\alpha_0=1\Biggr\}\quad
\text{if $k\geqq2$.}
\end{align*}
We give $G_p^{\langle k\rangle}(A_*)$ a group structure by the composition of formal power series 
if $p=2$ or $k\geqq2$.  
If $p$ is an odd prime, we define a group structure of $G_p^{\langle1\rangle}(A_*)$ as follows. 
Define a subset $\bar G_p^{\langle1\rangle}(A_*)$ of $A_*[\epsilon]/(\epsilon^2)[[X]]$ by 
$$\bar G_p^{\langle1\rangle}(A_*)=
\Biggl\{\sum\limits_{i=0}^{\infty}\alpha_iX^{p^i}\in A_*[\epsilon]/(\epsilon^2)[[X]]\,\Biggl|\,
\deg\,\alpha_i=2(p^{i+1}-p)\;(i\geqq0),\,\alpha_0-1\in(\epsilon)\Biggr\}.$$
We regard $A^*$ as a quotient of $A_*[\epsilon]/(\epsilon^2)$ with quotient map 
$q_{A_*}:A_*[\epsilon]/(\epsilon^2)\to A_*[\epsilon]/(\epsilon)=A_*$ and define a map 
$\hat q_{A_*}:\bar G_p^{\langle1\rangle}(A_*)\to G_p^{\langle1\rangle}(A_*)$ by 
$\hat q_{A_*}(\alpha(X))=\alpha_0X+\sum\limits_{i=1}^{\infty}q_{A_*}(\alpha_i)X^{p^i}$ if 
$\alpha(X)=\sum\limits_{i=0}^{\infty}\alpha_iX^{p^i}$. 
For $\alpha(X),\beta(X)\in G_p^{\langle1\rangle}(A_*)$, the product $\alpha(X)\cdot\beta(X)$ is 
defined to be $\hat q_{A_*}(\beta(\alpha(X)))$. 

Define subalgebras $\ca_p\langle k\rangle_*$ of $\ca_{p*}$ by 
$$\ca_p\langle k\rangle_*=\begin{cases}
\bff_2\bigl[\zeta_1^{2^k},\zeta_2^{2^k},\dots,\zeta_n^{2^k},\dots\bigr]& \text{if $p=2$}\\[2mm]
E(\tau_0)\otimes_{\bff_p}\!\bff_p\bigl[\xi_1^p,\xi_2^p,\dots,\xi_n^p,\dots\bigr]& 
\text{if $p$ is an odd prime and $k=1$}\\[2mm]
\bff_p\bigl[\xi_1^{p^k},\xi_2^{p^k},\dots,\xi_n^{p^k},\dots\bigr]& 
\text{if $p$ is an odd prime and $k\geqq2$.}
\end{cases}$$
Then $\ca_p\langle k\rangle_*$ is a Hopf subalgebra of $\ca_{p*}$. 
We put $\ca_p\langle0\rangle_*=\ca_{p*}$ and $G_p^{\langle0\rangle}=G_p$ for convenience. 

For a graded commutative algebra $R_*$, we denote by $h_{R_*}$ the functor represented by $R_*$ below.  

\begin{Prop}\label{G_p<k>}
$\ca_p\langle k\rangle_*$ represents $G_p^{\langle k\rangle}$. 
\end{Prop}

\begin{Proof}
Let $A_*$ be a graded commutative algebra over $\bff_p$ with product 
$m_{A_*}:A_*\otimes_{\bff_p}A_*\to A_*$. 
For $\varphi\in h_{\ca_p\langle k\rangle_*}(A_*)$, we define an element $\alpha_{\varphi}(X)$ of 
$G_p^{\langle k\rangle}(A_*)$ as follows. 
$$\alpha_{\varphi}(X)=\begin{cases}
X+\sum\limits_{i=1}^{\infty}\varphi\bigl(\zeta_i^{2^k}\bigr)X^{2^i}&\text{if $p=2$} \\ 
(1+\varphi(\tau_0)\epsilon)X+\sum\limits_{i=1}^{\infty}\varphi(\xi_i^p)X^{p^i}&
\text{if $p$ is an odd prime and $k=1$}\\ 
X+\sum\limits_{i=1}^{\infty}\varphi\bigl(\xi_i^{p^k}\bigr)X^{p^i}&
\text{if $p$ is an odd prime and $k\geqq2$}\end{cases}$$ 
Then, a correspondence $\varphi\mapsto\alpha_{\varphi}(X)$ gives a natural equivalence 
$\theta_k:h_{\ca_p\langle k\rangle_*}\to G_p^{\langle k\rangle}$. 

For $\varphi,\psi\in h_{\ca_p\langle k\rangle_*}(A_*)$, since the product $\psi\varphi$ of $\varphi$ 
is a composition 
$$\ca_p\langle k\rangle_*\xrightarrow{\mu}\ca_p\langle k\rangle_*\otimes_{\bff_p}
\ca_p\langle k\rangle_*\xrightarrow{\psi\otimes\varphi}A_*\otimes_{\bff_p}A_*\xrightarrow{m_{A_*}}
A_*,$$
it can be verified from (1) of (\ref{steenrod group rem}) and formulas of the coproduct of $\ca_{p*}$ 
given in (\cite{[17]}, See Intorduction above) that $\alpha_{\psi}(X)\cdot\alpha_{\varphi}(X)=
\alpha_{\psi\varphi}(X)$. 
Hence $\theta_k:h_{\ca_p\langle k\rangle_*}\to G_p^{\langle k\rangle}$ preserves group structures. 
\end{Proof}

We denote by $\iota_0:\ca_p\langle1\rangle_*\to\ca_{p*}$ and 
$\iota_k:\ca_p\langle k+1\rangle_*\to\ca_p\langle k\rangle_*$ ($k\geqq1$) the inclusion maps. 
It is well-known that inclusion morphisms between ungraded Hopf algebras over a field are faithfully 
flat (\cite{WCW}, Chapter 14). 
Since we are working in the category of graded algebras, we give a direct proof of the faithful 
flatness of the inclusion morphisms between graded Hopf algebras in the following special case. 

\begin{Prop}\label{faithful flatness of iota_k}
$\iota_0:\ca_p\langle1\rangle_*\to\ca_{p*}$ and 
$\iota_k:\ca_p\langle k+1\rangle_*\to\ca_p\langle k\rangle_*$ are faithfully flat for $k\geqq1$. 
\end{Prop}

\begin{Proof}
Let $V\langle k\rangle_{2,n}$ be a $2$-dimensional subspace of $\ca_2\langle k\rangle_*$ spanned by 
$1,\zeta_n^{2^k}$ and $V\langle k\rangle_{p,n}$ a $p$-dimensional subspace of 
$\ca_p\langle k\rangle_*$ spanned by $1,\xi_n^{p^k},\xi_n^{2p^k},\dots,\xi_n^{(p-1)p^k}$ if $p$ is 
an odd prime. 
Then, $\ca_p\langle k\rangle_*$ is isomorphic to 
$\biggl(\bigotimes\limits_{n\geqq1}V\langle k\rangle_{p,n}\biggr)\!\otimes_{\bff_p}
\ca_p\langle k+1\rangle_*$ for $k\geqq2$ as an $\ca_p\langle k+1\rangle_*$-module, 
$\ca_p\langle1\rangle_*$ is isomorphic to 
$E(\tau_0)\otimes_{\bff_p}\biggl(\bigotimes\limits_{n\geqq1}V\langle1\rangle_{p,n}\biggr)\!
\otimes_{\bff_p}\ca_p\langle2\rangle_*$ as an $\ca_p\langle2\rangle_*$-module and 
$\ca_{p*}$ is isomorphic to $E(\tau_1,\tau_2,\dots,\tau_n,\dots)\otimes_{\bff_p}\!
\biggl(\bigotimes\limits_{n\geqq1}V\langle0\rangle_{p,n}\biggr)\!\otimes_{\bff_p}
\ca_p\langle1\rangle_*$ as an $\ca_p\langle1\rangle_*$-module for an odd prime $p$. 
Hence $\ca_p\langle k\rangle_*$ is a free $\ca_p\langle k+1\rangle_*$-module and 
$\ca_{p*}$ is a free $\ca_p\langle 1\rangle_*$-module. 
It follows that the inclusion maps $\iota_0$ and $\iota_k$ are flat. 
Since $\ca_{2*}$ is integral over $\ca_2\langle1\rangle_*$ and 
$\ca_p\langle k\rangle_*$ is integral over $\ca_p\langle k+1\rangle_*$ if $k\geqq2$, 
$\spec(\iota_0):\spec\ca_{2*}\to\spec\ca_2\langle1\rangle_*$ and 
$\spec(\iota_k):\spec\ca_p\langle k\rangle_*\to\spec\ca_p\langle k+1\rangle_*$ are surjective for 
$k\geqq2$. 
Since $\tau_n$'s are nilpotent, the inclusion maps 
$E(\tau_0)\otimes_{\bff_p}\bff_p[\xi_1,\xi_2,\dots,\xi_n,\dots]\to\ca_{p*}$ and 
$\bff_p[\xi_1^p,\xi_2^p,\dots,\xi_n^p,\dots]\to\ca_p\langle1\rangle_*$ induces bijections 
$\spec\ca_{p*}\to\spec(E(\tau_0)\otimes_{\bff_p}\!\bff_p[\xi_1,\xi_2,\dots,\xi_n,\dots])$ and 
$\spec\ca_p\langle1\rangle_*\to\spec\bff_p[\xi_1^p,\xi_2^p,\dots,\xi_n^p,\dots]$, respectively.  
Since $E(\tau_0)\otimes_{\bff_p}\!\bff_p[\xi_1,\xi_2,\dots,\xi_n,\dots]$ is integral over 
$\ca_p\langle1\rangle_*$ and $\bff_p[\xi_1^p,\xi_2^p,\dots,\xi_n^p,\dots]$ is integral over 
$\ca_p\langle2\rangle_*$, the following maps induced by the inclusion maps are surjective. 
$$\spec(E(\tau_0)\otimes_{\bff_p}\!\bff_p[\xi_1,\xi_2,\dots,\xi_n,\dots])\to 
\spec\ca_p\langle1\rangle_*, \quad
\spec\bff_p[\xi_1^p,\xi_2^p,\dots,\xi_n^p,\dots]\to \spec\ca_p\langle2\rangle_*$$ 
Therefore $\spec(\iota_0):\spec\ca_{p*}\to\spec\ca_p\langle1\rangle_*$ and 
$\spec(\iota_1):\spec\ca_p\langle1\rangle_*\to\spec\ca_p\langle2\rangle_*$ are surjective. 
It follows from Theorem 3 of \cite{ca} that $\iota_k$'s are faithfully flat.  
\end{Proof}

We define morphisms $\rho_0:G_p\to G_p^{\langle1\rangle}$ and 
$\rho_k:G_p^{\langle k\rangle}\to G_p^{\langle k+1\rangle}$ of affine group schemes by 
$\rho_{0A_*}(\alpha(X))=\alpha_0X+\sum\limits_{i=1}^{\infty}\alpha_i^pX^{p^i}$ and 
$\rho_{kA_*}(\alpha(X))=X+\sum\limits_{i=1}^{\infty}\alpha_i^pX^{p^i}$, respectively if 
$\alpha(X)=\sum\limits_{i=0}^{\infty}\alpha_iX^{p^i}$. 
We consider the natural equivalence $\theta_k:h_{\ca_p\langle k\rangle_*}\to G_p^{\langle k\rangle}$  
given in the proof of (\ref{G_p<k>}) and the natural equivalence $\theta:h_{\ca_{p*}}\to G_p$ given 
in Proposition 5.1 of \cite{ay1}. 

\begin{Prop}\label{abelian kernel1}
Let us denote by $\iota_0^*:h_{\ca_{p*}}\to h_{\ca_p\langle1\rangle_*}$ and 
$\iota_k^*:h_{\ca_p\langle k\rangle_*}\to h_{\ca_p\langle k+1\rangle_*}$ the natural transformations 
induced by $\iota_0$ and $\iota_k$, respectively. 
The following diagrams are commutative. 
$$\begin{CD}
h_{\ca_{p*}} @>{\iota_0^*}>> h_{\ca_p\langle1\rangle_*} \\
@VV{\theta}V @VV{\theta_1}V \\
G_p @>{\rho_0}>> G_p^{\langle1\rangle}
\end{CD}
\qquad
\begin{CD}
h_{\ca_p\langle k\rangle_*} @>{\iota_k^*}>> h_{\ca_p\langle k+1\rangle_*} \\
@VV{\theta_k}V @VV{\theta_{k+1}}V \\
G_p^{\langle k\rangle} @>{\rho_{k}}>> G_p^{\langle k+1\rangle}
\end{CD}$$
\end{Prop}

\begin{Proof}
For a graded commutative algebra $A_*$ and $\varphi\in h_{\ca_p\langle k\rangle_*}(A_*)$, we have an 
equality 
$$\rho_{kA_*}(\alpha_{\varphi}(X))=
\begin{cases}X+\sum\limits_{i=1}^{\infty}\varphi\bigl(\zeta_i^{2^{k+1}}\bigr)X^{2^i}&\text{if $p=2$}
\\ 
X+\sum\limits_{i=1}^{\infty}\varphi\bigl(\xi_i^{p^{k+1}}\bigr)X^{p^i}&
\text{if $p$ is an odd prime}\end{cases}$$
which shows $\rho_{kA_*}(\theta_{kA_*}(\varphi))=\rho_{kA_*}(\alpha_{\varphi}(X))=
\alpha_{\varphi\iota_k}(X)=\theta_{k+1A_*}(\varphi\iota_k)$. 

Similarly, since $\theta_{A_*}:h_{\ca_{p*}}(A_*)\to G_p(A_*)$ is given by 
$$\theta_{A_*}(\varphi)=\begin{cases}
X+\sum\limits_{i=1}^{\infty}\varphi(\zeta_i)X^{2^i}&\text{if $p=2$} \\ 
(1+\varphi(\tau_0)\epsilon)X+
\sum\limits_{i=1}^{\infty}(\varphi(\xi_i)+\varphi(\tau_i)\epsilon)X^{p^i}
&\text{if $p$ is an odd prime,}\end{cases}$$
we have the following equality.  
$$\rho_{0A_*}(\theta_{A_*}(\varphi))=\begin{cases}
X+\sum\limits_{i=1}^{\infty}\varphi(\zeta_i^2)X^{2^i}&\text{if $p=2$} \\ 
(1+\varphi(\tau_0)\epsilon)X+\sum\limits_{i=1}^{\infty}\varphi(\xi_i^p)X^{p^i}
&\text{if $p$ is an odd prime}\end{cases}$$ 
Hence $\rho_{0A_*}(\theta_{A_*}(\varphi))=\theta_{1A_*}(\varphi\iota_0)$ holds for 
$\varphi\in h_{\ca_{p*}}(A_*)$. 
\end{Proof}

Let us denote by $G_p^{ab}$ the kernel of $\rho_0:G_p\to G_p^{\langle1\rangle}$ and by $G_p^{[k]}$ 
the kernel of $\rho_k:G_p^{\langle k\rangle}\to G_p^{\langle k+1\rangle}$. 
The following fact is a direct consequence of the definitions of $\rho_0$ and $\rho_k$. 

\begin{Prop}\label{abelian kernel2}
Let $A_*$ be a graded commutative $\bff_p$-algebra. 

$(1)$ $G_p^{ab}(A_*)$ is identified with a subgroup of $G_p(A_*)$ given by 
$$\biggl\{\alpha(X)\in G_p(A_*)\,\bigg|\,\alpha(X)=X+\sum\limits_{i=1}^{\infty}\alpha_iX^{p^i},\;
\alpha_i^p=0\; \text{for}\;i\geqq1\biggr\}.$$

$(2)$ $G_p^{[1]}(A_*)$ is identified with a subgroup of $G_p^{\langle1\rangle}(A_*)$ given by 
$$\biggl\{\alpha(X)\in G_p^{\langle1\rangle}(A_*)\,\bigg|\,\alpha(X)=
\alpha_0X+\sum\limits_{i=1}^{\infty}\alpha_iX^{p^i},\;\alpha_i^p=0\; \text{for}\;i\geqq1\biggr\}.$$

$(3)$ $G_p^{[k]}(A_*)$ $(k\geqq2)$ is identified with a subgroup of $G_p^{\langle k\rangle}(A_*)$ 
given by 
$$\biggl\{\alpha(X)\in G_p^{\langle k\rangle}(A_*)\,\bigg|\,\alpha(X)=
X+\sum\limits_{i=1}^{\infty}\alpha_iX^{p^i},\;\alpha_i^p=0\; \text{for}\;i\geqq1\biggr\}.$$
\end{Prop}

Let $I_p^{\langle k\rangle}$ be an ideal of $\ca_p\langle k\rangle_*$ ($\ca_{p*}$ if $k=0$) given by 
$$I_p^{\langle k\rangle}=\begin{cases}
\bigl(\zeta_1^{2^{k+1}},\zeta_2^{2^{k+1}},\cdots,\zeta_n^{2^{k+1}},\cdots\bigr) 
&\text{if $p=2$}\\[2mm]
\bigl(\tau_0,\xi_1^p,\xi_2^p,\cdots,\xi_n^p,\cdots\bigr)&
\text{if $p$ is an odd prime and $k=0$}\\[2mm]
\bigl(\xi_1^{p^{k+1}},\xi_2^{p^{k+1}},\cdots,\xi_n^{p^{k+1}},\cdots\bigr)&
\text{if $p$ is an odd prime and $k\geqq1$.}
\end{cases}$$ 

\begin{Prop}\label{abelian kernel}
$G_p^{ab}$ is an abelian group scheme represented by $\ca_{p*}/I_p^{\langle0\rangle}$. 
Similarly, $G_p^{[k]}$ is an abelian group scheme represented by 
$\ca_p\langle k\rangle_*/I_p^{\langle k\rangle}$. 
\end{Prop}

\begin{Proof}
It follows from (\ref{abelian kernel1}) that $G_p^{ab}$ is represented by 
$\ca_{p*}/I_p^{\langle0\rangle}$ and that $G_p^{[k]}$ is represented by 
$\ca_p\langle k\rangle_*/I_p^{\langle k\rangle}$. 
The formulas of the coproduct of $\ca_{p*}$ and the following equalities imply that 
$\ca_{p*}/I_p^{\langle0\rangle}$ and $\ca_p\langle k\rangle_*/I_p^{\langle k\rangle}$ are 
primitively generared. 
Hence $G_p^{ab}$ and $G_p^{[k]}$ are abelian group schemes. 
\end{Proof}

\begin{Remark}\label{abelian kernel rem}
We denote by $E(n)_*$ a monogenic Hopf graded algebra over $\bff_p$ generated by 
$\tau_n$ with $\tau_n^2=0$ and $\deg\,\tau_n=2p^n-1$. 
We also denote by $F(n,k)_*$ a monogenic Hopf graded algebra over $\bff_p$ generated by 
$\xi_{n,k}$ with $\xi_{n,k}^p=0$ and $\deg\,\xi_{n,k}=2p^k(p^n-1)$. 
Then, $\ca_{p*}/I_p^{\langle0\rangle}$ is isomorphic to $\biggl(\bigotimes\limits_{n\geqq1}
E(n)_*\biggr)\otimes_{\bff_p}\!\biggl(\bigotimes\limits_{n\geqq1}F(n,0)_*\biggr)$, 
$\ca_p\langle1\rangle_*/I_p^{\langle1\rangle}$ is isomorphic to 
$E(0)_*\otimes_{\bff_p}\!\biggl(\bigotimes\limits_{n\geqq1}F(n,1)_*\biggr)$ if $p$ is an odd prime 
and $\ca_p\langle k\rangle_*/I_p^{\langle k\rangle}$ is isomorphic to 
$\bigotimes\limits_{n\geqq1}F(n,k)_*$ if $p=2$ or $k\geqq2$. 
\end{Remark}

For a graded vector space $M_*$ over $\bff_p$, we denote by 
$T_{M_*}:M_*\otimes_{\bff_p}M_*\to M_*\otimes_{\bff_p}M_*$ the switching map given by 
$T_{M_*}(x\otimes y)=(-1)^{\deg\,x\,\deg\,y}y\otimes x$. 

\begin{Prop}\label{maximal abelian subgroup}
If $I$ is a Hopf ideal of $\ca_{p*}$ contained in $I_p^{\langle0\rangle}$ and $\ca_{p*}/I$ is a 
cocomutative Hopf algebra, then $I=I_p^{\langle0\rangle}$ holds. 
Similarly, if $I$ is a Hopf ideal of $\ca_p\langle k\rangle_*$ contained in $I_p^{\langle k\rangle}$ 
and $\ca_p\langle k\rangle_*/I$ is a cocomutative Hopf algebra, then $I=I_p^{\langle k\rangle}$ 
holds. 
\end{Prop}

\begin{Proof}
Let us denote by $\pi_I:\ca_{p*}\to\ca_{p*}/I$ the quotient map and by 
$\mu_I:\ca_{p*}/I\to\ca_{p*}/I\otimes_{\bff_p}\ca_{p*}/I$ the coproduct of $\ca_{p*}/I$ induced by 
the coproduct $\mu$ of $\ca_{p*}$. 
Since $\mu_I=T_{\ca_{p*}/I}\mu_I$ by the assumption, the image of 
$\mu-T_{\ca_{p*}}\mu:\ca_{p*}\to\ca_{p*}\otimes_{\bff_p}\ca_{p*}$ is contained in the kernel of 
$\pi_I\otimes_{\bff_p}\pi_I$ which is $\ca_{p*}\otimes_{\bff_p}I+I\otimes_{\bff_p}\ca_{p*}$. 
$$\begin{CD}
\ca_{p*} @>{\phantom{MM}\mu-T_{\ca_{p*}}\mu\phantom{MM}}>> \ca_{p*}\otimes_{\bff_p}\ca_{p*} \\
@VV{\pi_I}V @VV{\pi_I\otimes_{\bff_p}\pi_I}V \\
\ca_{p*}/I @>{\mu_I-T_{\ca_{p*}/I}\mu_I\,=\,0}>> \ca_{p*}/I\otimes_{\bff_p}\ca_{p*}/I
\end{CD}$$
Consider the case that $p$ is an odd prime and assume $\tau_0\not\in I$. 
Then, since $I\subset I_p^{\langle0\rangle}$, $I$ does not contain any element whose degree is less 
than $2p(p-1)$. 
It follows that $\pi_I$ is injective in degree less than $2p(p-1)$, hence so is 
$\pi_I\otimes_{\bff_p}\pi_I$. 
However, a non-zero element $\xi_1\otimes\tau_0-\tau_0\otimes\xi_1=
\mu(\tau_1)-T_{\ca_{p*}}\mu(\tau_1)$ of $\ca_{p*}\otimes_{\bff_p}\ca_{p*}$ has degree $2p-1$ 
but it is contained in the kernel of $\pi_I\otimes_{\bff_p}\pi_I$, which is a contradiction. 
Therefore $I$ containes $\tau_0$. 
Assume inductively that $\xi_1^p,\xi_2^p,\dots,\xi_{i-1}^p\in I$ for $i\geqq1$. 
The following relations imply   
$\xi_i^p\otimes\xi_1-\xi_1\otimes\xi_i^p\in I\otimes_{\bff_p}\ca_{p*}+\ca_{p*}\otimes_{\bff_p}I$. 
\begin{align*}
\sum\limits_{j=1}^i\bigl(\xi_{i+1-j}^{p^j}\otimes\xi_j-\xi_j\otimes\xi_{i+1-j}^{p^j}\bigr)=
\mu(\xi_{i+1})-T_{\ca_{p*}}\mu_I(\xi_{i+1})&\in I\otimes_{\bff_p}\ca_{p*}+\ca_{p*}\otimes_{\bff_p}I
\\
\sum\limits_{j=2}^i\bigl(\xi_{i+1-j}^{p^j}\otimes\xi_j-\xi_j\otimes\xi_{i+1-j}^{p^j}\bigr)&\in 
I\otimes_{\bff_p}\ca_{p*}+\ca_{p*}\otimes_{\bff_p}I
\end{align*} 
Since $I\subset I_p^{\langle0\rangle}$, $I$ does not have any element of degree $2p-2$. 
Thus we have $\xi_i^p\in I$. 
Similarly, we can show that $\zeta_i^2\in I$ by the induction on $i$ if $p=2$. 
 
Let us denote by $\pi^k_I:\ca_p\langle k\rangle_*\to\ca_p\langle k\rangle_*/I$ and $\mu^k_I:
\ca_p\langle k\rangle_*/I\to\ca_p\langle k\rangle_*/I\otimes_{\bff_p}\ca_p\langle k\rangle_*/I$
the quotient map and the coproduct of $\ca_p\langle k\rangle_*/I$ induced by the coproduct $\mu$ of 
$\ca_p\langle k\rangle_*$, respectively. 
Since $\mu^k_I=T_{\ca_p\langle k\rangle_*/I}\mu^k_I$ by the assumption, the image of 
$\mu-T_{\ca_p\langle k\rangle_*}\mu:\ca_p\langle k\rangle_*\to\ca_p\langle k\rangle_*\otimes_{\bff_p}
\ca_p\langle k\rangle_*$ is contained in the kernel of $\pi_I\otimes_{\bff_p}\pi_I$ 
which is $\ca_p\langle k\rangle_*\otimes_{\bff_p}I+I\otimes_{\bff_p}\ca_p\langle k\rangle_*$. 
$$\begin{CD}
\ca_p\langle k\rangle_* @>{\phantom{MM}\mu-T_{\ca_p\langle k\rangle_*}\mu\phantom{MM}}>> 
\ca_p\langle k\rangle_*\otimes_{\bff_p}\ca_p\langle k\rangle_* \\
@VV{\pi_I}V @VV{\pi_I\otimes_{\bff_p}\pi_I}V \\
\ca_p\langle k\rangle_*/I @>{\mu_I-T_{\ca_p\langle k\rangle_*/I}\mu_I\,=\,0}>> 
\ca_p\langle k\rangle_*/I\otimes_{\bff_p}\ca_p\langle k\rangle_*/I
\end{CD}$$
Consider the case that $p$ is an odd prime. 
Assume inductively that $\xi_1^{p^{k+1}},\xi_2^{p^{k+1}},\dots,\xi_{i-1}^{p^{k+1}}\in I$ for 
$i\geqq1$. 
The following relations imply   
$\xi_i^{p^{k+1}}\otimes\xi_1^{p^k}-\xi_1^{p^k}\otimes\xi_i^{p^{k+1}}\in 
I\otimes_{\bff_p}\ca_p\langle k\rangle_*+\ca_p\langle k\rangle_*\otimes_{\bff_p}I$. 
\begin{align*}
\sum\limits_{j=1}^i
\bigl(\xi_{i+1-j}^{p^{j+k}}\otimes\xi_j^{p^k}-\xi_j^{p^k}\otimes\xi_{i+1-j}^{p^{j+k}}\bigr)=
\mu\bigl(\xi_{i+1}^{p^k}\bigr)-T_{\ca_p\langle j\rangle_*/I}\mu_I\bigl(\xi_{i+1}^{p^k}\bigr)
&\in I\otimes_{\bff_p}\ca_p\langle j\rangle_*+\ca_p\langle j\rangle_*\otimes_{\bff_p}I
\\
\sum\limits_{j=2}^i
\bigl(\xi_{i+1-j}^{p^{j+k}}\otimes\xi_j^{p^k}-\xi_j^{p^k}\otimes\xi_{i+1-j}^{p^{j+k}}\bigr)
&\in I\otimes_{\bff_p}\ca_p\langle j\rangle_*+\ca_p\langle j\rangle_*\otimes_{\bff_p}I
\end{align*} 
Since $I\subset I_p^{\langle k\rangle}$, $I$ does not have any element of degree less than 
$2p^{k+1}(p-1)$. 
Thus we have $\xi_i^{p^{k+1}}\in I$. 
Similarly, we can show that $\zeta_i^{2^{k+1}}\in I$ by the induction on $i$ if $p=2$. 
\end{Proof}

The above result shows that $I_p^{\langle0\rangle}$ is a minimal Hopf ideal of $\ca_{p*}$ such that 
the quotient Hopf algebra $\ca_{p*}/I_p^{\langle0\rangle}$ is cocomutative, hence it follows from 
(\ref{steenrod group2}) that $G_p^{ab}$ is a maximal abelian subgroup of $G_p$. 
Similarly, $I_p^{\langle k\rangle}$ is a minimal Hopf ideal of $\ca_p\langle k\rangle_*$ such that 
the quotient Hopf algebra $\ca_p\langle k\rangle_*/I_p^{\langle0\rangle}$ is cocomutative, thus 
$G_p^{[k]}$ is a maximal abelian subgroup of $G_p^{\langle k\rangle}$ by (\ref{G_p<k>}). 

Summarizing the results in this section so far, we have the following. 

\begin{Thm}\label{tower}
$(1)$ There is a short exact sequence $1\to G_p^{ab}\xrightarrow{\kappa_0}G_p\xrightarrow{\rho_0}
G_p^{\langle1\rangle}\to1$ of affine group schemes over $\bff_p$ such that $G_p^{ab}$ is a maximal 
abelian subgroup scheme of $G_p$ and $\rho_0$ is faithfully flat. 

$(2)$ There is a short exact sequence $1\to G_p^{[k]}\xrightarrow{\kappa_k}G_p^{\langle k\rangle}
\xrightarrow{\rho_k}G_p^{\langle k+1\rangle}\to1$ of affine group schemes over $\bff_p$ such that 
$G_p^{[k]}$ is a maximal abelian subgroup scheme of $G_p^{\langle k\rangle}$ and $\rho_k$ is 
faithfully flat. 
$$\begin{CD}
G_p^{ab} @. G_p^{[1]} @. G_p^{[2]} @.\cdots @. G_p^{[k]} @. G_p^{[k+1]} @. \cdots \\
@VV{\kappa_0}V @VV{\kappa_1}V @VV{\kappa_2}V @. @VV{\kappa_k}V @VV{\kappa_{k+1}}V \\
G_p @>{\rho_0}>> G_p^{\langle1\rangle} @>{\rho_1}>> G_p^{\langle2\rangle} @>{\rho_2}>> \cdots 
@>{\rho_{k-1}}>> G_p^{\langle k\rangle} @>{\rho_k}>> G_p^{\langle k+1\rangle} @>{\rho_{k+1}}>> \cdots
\end{CD}$$
\end{Thm}

We denote by $J_p^{\langle k\rangle}$ the ideal of $\ca_{p*}$ that is the extension of 
$I_p^{\langle k\rangle}$ by the inclusion map $\ca_p\langle k\rangle_*\to\ca_{p*}$. 
Then, we have a decreasing sequence 
$I_p^{\langle0\rangle}=J_p^{\langle0\rangle}\supset J_p^{\langle1\rangle}\supset\cdots\supset 
J_p^{\langle k\rangle}\supset J_p^{\langle k+1\rangle}\supset\cdots$ of Hopf ideals of $\ca_{p*}$. 
Hence if we denote by $\tilde G_p^{\langle k\rangle}$ the closed subgroup scheme of $G_p$ defined by 
$J_p^{\langle k\rangle}$, we have an increasing sequence $G_p^{ab}=\tilde G_p^{\langle0\rangle}
\subset\tilde G_p^{\langle1\rangle}\subset\cdots\subset\tilde G_p^{\langle k\rangle}\subset
\tilde G_p^{\langle k+1\rangle}\subset\cdots\subset G_p$ of subgroup schemes of $G_p$. 
In fact, $\tilde G_p^{\langle k\rangle}$ is the inverse image of $G_p^{[k]}$ by a composition 
$G_p\xrightarrow{\rho_0}G_p^{\langle1\rangle}\xrightarrow{\rho_1}G_p^{\langle2\rangle}
\xrightarrow{\rho_2}\cdots\xrightarrow{\rho_{k-1}}G_p^{\langle k\rangle}$, which implies that 
$\tilde G_p^{\langle k\rangle}$ is a normal subgroup scheme of $G_p$. 
Moreover, the map $\ca_p\langle k\rangle_*/I_p^{\langle k\rangle}\to\ca_{p*}/J_p^{\langle k\rangle}$ 
induced by the inclusion map $\ca_p\langle k\rangle_*\to\ca_{p*}$ is faithfully flat whose images 
of elements of positive degree generate an ideal $J_p^{\langle k-1\rangle}/J_p^{\langle k\rangle}$ of 
$\ca_{p*}/J_p^{\langle k\rangle}$. 
Thus we see the following fact. 

\begin{Prop}\label{another filtration}
There is an increasing sequence 
$$G_p^{ab}=\tilde G_p^{\langle0\rangle}\subset\tilde G_p^{\langle1\rangle}\subset\cdots\subset
\tilde G_p^{\langle k\rangle}\subset\tilde G_p^{\langle k+1\rangle}\subset\cdots\subset G_p$$ 
of normal subgroup schemes of $G_p$ such that 
$\tilde G_p^{\langle k\rangle}/\tilde G_p^{\langle k-1\rangle}$ is isomorphic to $G_p^{[k]}$. 
\end{Prop}

Finally, we investigate the dual of the successive quotients 
$$\begin{CD}
\cdots @= \ca_{p*} @= \ca_{p*} @= \cdots @= \ca_{p*} @= \ca_{p*} \\
@. @VVV @VVV @. @VVV @VVV \\
\cdots @>>> \ca_{p*}/J_p^{\langle k+1\rangle} @>>> \ca_{p*}/J_p^{\langle k\rangle} @>>> \cdots @>>> 
\ca_{p*}/J_p^{\langle1\rangle} @>>> \ca_{p*}/J_p^{\langle0\rangle}
\end{CD}$$
of the dual Steenrod algebra. 
To do this, we recall the Milnor basis given in \cite{[17]} briefly. 

Let $\seq$ be the set of all infinite sequences  $(r_1,r_2,\dots,r_i,\dots)$ of 
non-negative integers such that $r_i=0$ for all but finite number of $i$'s. 
$\seq$ is regarded as an abelian monoid with unit $\vzr=(0,0,\dots)$ by componentwise addition. 
We denote by $\seq^b$ a subset of $\seq$ consisting of all sequences $(e_0,e_1,\dots,e_i,\dots)$ 
such that $e_i=0,1$ for all $i=0,1,\dots$. 
Define an order $\leqq$ in $\seq$ by ``$R\leqq S$ if and only if $R=(r_1,r_2,\dots,r_i,\dots),
S=(s_1,s_2,\dots,s_i,\dots)\in\seq$ and $r_i\leqq s_i$ for all $i=1,2,\dots$''. 
If $r_i=0$ for $i>n$, we denote $(r_1,r_2,\dots,r_i,\dots)$ by $(r_1,r_2,\dots,r_n)$. 
For a positive integer $n$, let us denote by $E_n$ an element of $\seq^b$ such that the $n$-th 
entry is $1$ and other entries are all zero. 

For $E=(e_0,e_1,\dots,e_m)\in\seq^b$ and $R=(r_1,r_2,\dots,r_n)\in\seq$, we put 
$\tau(E)=\tau_0^{e_0}\tau_1^{e_1}\cdots\tau_m^{e_m}$,  
$\xi(R)=\xi_1^{r_1}\xi_2^{r_2}\cdots\xi_n^{r_n}$, 
$\zeta(R)=\zeta_1^{r_1}\zeta_2^{r_2}\cdots\zeta_i^{r_n}$. 
We consider the monomial basis $B_2=\{\zeta(R)\,|\,R\in\seq\}$ of $\ca_{2*}$ and 
$B_p=\{\tau(E)\xi(R)\,|\,E\in\seq^b,\;R\in\seq\}$ of $\ca_{p*}$ if $p$ is an odd prime. 
The following assertion is clear from the definition of $J_p^{\langle k\rangle}$. 

\begin{Lemma}\label{basis of J_p}
$\bigl\{\zeta(R)\,\big|\,R\geqq 2^{k+1}E_n\;\;\text{for some}\;\;n\geqq1\bigr\}$ is a basis of 
$J_2^{\langle k\rangle}$. 
If $p$ is an odd prime, then $\{\tau(E)\xi(R)\,|\,E\in\seq^b\;\text{and}\;\;
``E\geqq E_1\;\text{or}\;\;R\geqq pE_n\;\;\text{for some}\;\;n\geqq1"\bigr\}$ 
is a basis of $J_p^{\langle0\rangle}$ and, for $k\geqq1$,  
$\{\tau(E)\xi(R)\,|\,E\in\seq^b\;\text{and}\;\;R\geqq p^{k+1}E_n\;\;\text{for some}\;\;n\geqq1
\bigr\}$ is a basis of $J_p^{\langle k\rangle}$ . 
\end{Lemma}

Let $\sq(R)$ be the dual of $\zeta(R)$ with respect to $B_2$ if $p=2$. 
If $p$ is an odd prime, let us denote by $\wp(R)$ and $Q_i$ the duals of $\xi(R)$ and $\tau_i$ with 
respect to $B_p$, respectively. 
We put $Q(E)=Q_0^{e_0}Q_1^{e_1}\cdots Q_m^{e_m}$ if $E=(e_0,e_1,\dots,e_m)\in\seq^b$. 
Let us denote by $\langle\theta,x\rangle$ the Kronecker pairing of $\theta\in\ca_p^*$ and 
$x\in\ca_{p*}$. 
Since the dual of $\ca_{p*}/J_p^{\langle k\rangle}$ is canonically isomorphic to a subspace 
$\{\theta\in\ca_p^*\,|\,\langle\theta,x\rangle=0\;\text{for all}\;x\in J_p^{\langle k\rangle}\}$ 
of $\ca_p^*$, it follows from (\ref{maximal abelian subgroup}), (\ref{basis of J_p}) and Theorem 4a 
of \cite{[17]} that we have the following result. 

\begin{Prop}\label{filtration of A_p^*}
Let $\ca_2\langle k\rangle^*$ be a subspace of $\ca_2^*$ spanned by the following set. 
$$\bigl\{\sq(r_1,r_2,\dots,r_i,\dots)\,\big|\,r_n<2^{k+1}\;\text{for all}\;\;n\geqq1\bigr\}$$
If $p$ is an odd prime, let $\ca_p\langle k\rangle^*$ be a subspace of $\ca_p^*$ spanned by the 
following set. 
\begin{align*}
&\bigl\{Q(e_0,e_1,\dots,e_i,\dots)\wp(r_1,r_2,\dots,r_i,\dots)\,\big|\,e_0=0\;\;\text{and}\;\;
r_n<p\;\;\text{for all}\;\;n\geqq1\bigr\} & \text{if $k=0$}
\\
&\bigl\{Q(e_0,e_1,\dots,e_i,\dots)\wp(r_1,r_2,\dots,r_i,\dots)\,\big|\,r_n<p^{k+1}\;
\text{for all}\;\;n\geqq1\bigr\} & \text{if $k\geqq1$}
\end{align*}
Here we assume that $(e_0,e_1,\dots,e_i,\dots)\in\seq^b$. 
Then, $\ca_p\langle k\rangle^*$ is a Hopf subalgebra of $\ca_p^*$ which is the dual of 
$\ca_{p*}/J_p^{\langle k\rangle}$. 
Thus we have an increasing filtration 
$$\ca_p\langle0\rangle^*\subset\ca_p\langle1\rangle^*\subset\cdots\subset
\ca_p\langle k\rangle^*\subset\ca_p\langle k+1\rangle^*\subset\cdots\subset\ca_p^*$$
of Hopf subalgebras of $\ca_p^*$ such that $\ca_p\langle0\rangle^*$ is a maximal commutative Hopf 
subalgebra of $\ca_p^*$. 
\end{Prop}

\end{document}